\documentclass[a4paper]{amsart}
\usepackage{amssymb}
\usepackage{amscd}
\usepackage{amsthm}
\usepackage{amsmath}
\usepackage{mathtools}
\usepackage{latexsym}
\usepackage{hyperref}
\usepackage[all]{xy}
\usepackage[utf8]{inputenc}
\usepackage{enumitem}

\usepackage{tikz-cd}

\setlist[enumerate]{label={(\roman*)}}

\theoremstyle{plain}
\newtheorem{theorem}{Theorem}
\newtheorem{corollary}[theorem]{Corollary}
\newtheorem{lemma}[theorem]{Lemma}
\newtheorem{proposition}[theorem]{Proposition}

\theoremstyle{definition}
\newtheorem{definition}[theorem]{Definition}
\newtheorem{example}[theorem]{Example}

\theoremstyle{remark}

\newtheorem{remark}{Remark}
\newtheorem{acknowledgment}[theorem]{Acknowledgment}
\numberwithin{theorem}{section}
\usepackage{newtxtext}
\usepackage{newtxmath}

\DeclareMathAlphabet\urwscr{U}{urwchancal}{m}{n}%
\DeclareMathAlphabet\rsfscr{U}{rsfso}{m}{n}
\DeclareMathAlphabet\euscr{U}{eus}{m}{n}
\DeclareFontEncoding{LS2}{}{}
\DeclareFontSubstitution{LS2}{stix}{m}{n}
\DeclareMathAlphabet\stixcal{LS2}{stixcal}{m} {n}
\newcommand{\Filt}[1]{\mbox{\rm Filt}(#1)}
\newcommand{\dFilt}[1]{\mbox{\rm dFilt}(#1)}

\newcommand{\Spec}[1]{\operatorname{Spec}(#1)}

\newcommand{\card}{\mbox{\rm{card\,}}}
\newcommand{\gen}{\mbox{\rm{gen\,}}}

\newcommand{\Add}{\mbox{\rm{Add\,}}}

\newcommand{\im}{\mbox{\rm{Im\,}}}

\newcommand{\Hom}[3]{\operatorname{Hom}_{#1}(#2,#3)}
\newcommand{\Ext}[4]{\operatorname{Ext}^{#1}_{#2}(#3,#4)}

\newcommand{\rmod}[1]{\mbox{\rm{Mod}--}{#1}}

\begin{document}

\title{Categoricity for transfinite extensions of modules}

\author{Jan Trlifaj}
\address{Charles University, Faculty of Mathematics and Physics, Department of Algebra \\
Sokolovsk\'{a} 83, 186 75 Prague 8, Czech Republic}
\email{trlifaj@karlin.mff.cuni.cz}

\begin{abstract} For each deconstructible class of modules $\mathcal D$, we prove that the categoricity of $\mathcal D$ in a big cardinal is equivalent to its categoricity in a tail of cardinals. We also prove Shelah's Categoricity Conjecture for $(\mathcal D, \preceq )$, where $(\mathcal D, \preceq )$ is any abstract elementary class of roots of Ext.
\end{abstract}

\date{\today}

\thanks{Research supported by GA\v CR 20-13778S, and by the Fields Institute within the Thematic Program on Set Theoretic Methods in Algebra, Dynamics and Geometry.}
	


\maketitle


The two basic concepts studied in this paper are the deconstructibility and the categoricity of classes of modules. 

Many classes of modules appearing in homological algebra are known to be deconstructible. This means that they can be obtained by transfinite extensions from a set of their elements. 
Each deconstructible class provides precovers (right approximations) of modules, and hence it makes it possible to develop relative homological algebra in the spirit of \cite[Chap.\ 8]{EJ}. 

In contrast, categoricity of such classes is a very strong property that is quite rare. So it is reasonable to expect that once categoricity does occur in a sufficiently big cardinal, it occurs in all sufficiently big cardinals, that is, in a tail of cardinals. Our main goal here is to verify this expectation for all deconstructible classes of modules. In the particular case of the AEC's of roots of Ext, $(\mathcal D, \preceq )$ in the sense of \cite{BET} and \cite[\S 10.3]{GT}, we show that {\lq\lq}sufficiently big{\rq\rq} stands for {\lq\lq}of cardinality $\geq \rho{^+}${\rq\rq}, where $\rho$ is the L\"owenheim-Skolem number of $(\mathcal D, \preceq)$, which is the least infinite cardinal $\kappa \geq \card R$ such that $\mathcal D$ is $\kappa^+$-deconstructible. Thus, Shelah's Categoricity Conjecture holds for all AEC's of roots of Ext.       

\subsection{Transfinite extensions and deconstructible classes}
In what follows, $R$ denotes an (associative, unital) ring and $\rmod R$ the category of all (unitary right $R$-) modules.

Given a class $\mathcal C \subseteq \rmod R$, $\Filt{\mathcal C}$ denotes the class of all \emph{transfinite extensions of modules from $\mathcal C$} (also called \emph{$\mathcal C$-filtered} modules). These are the modules $M$ that possess a \emph{$\mathcal C$-filtration}, i.e., an increasing chain of submodules, $\mathcal M = ( M_\alpha \mid \alpha \leq \sigma )$, such that $M_0 = 0$, $M_\sigma = M$, $M_\alpha = \bigcup_{\beta < \alpha} M_\beta$ for each limit ordinal $\alpha \leq \sigma$, and $M_{\alpha + 1}/M_{\alpha}$ is isomorphic to an element of $\mathcal C$ for each $\alpha < \sigma$. The ordinal $\sigma$ is called the \emph{length} of $\mathcal M$.

Let $\mathcal D$ be a class of modules. If $\kappa$ is an infinite cardinal such that $\mathcal D$ contains a subset $\mathcal C$ consisting of $< \kappa$-presented modules and $\mathcal D = \Filt{\mathcal C}$, then $\mathcal D$ is called \emph{$\kappa$-deconstructible}. $\mathcal D$ is \emph{deconstructible} provided it is $\kappa$-deconstructible for some infinite cardinal $\kappa$.  For example, if $\mathcal C = \{ R \}$, then $\Filt{\mathcal C}$ is the $\aleph_0$-deconstructible class of all free modules. Similarly, if $\kappa = \card R + \aleph_0$ and $\mathcal C$ denotes the class of all flat modules of cardinality $\leq \kappa$, then $\Filt{\mathcal C} = \mathcal F _0$ is the $\kappa ^+$-deconstructible class of all flat modules (cf.\ \cite[6.17]{GT}).  

Let $\mathcal C$ be any class of modules. Then the class $\Filt{\mathcal C}$ is closed under arbitrary direct sums and (transfinite) extensions, but it need not be closed under direct summands. We will denote by $\dFilt{\mathcal C}$ the class of all direct summands of modules in $\Filt{\mathcal C}$. So for example, $\dFilt{\mathcal C} = \mathcal P_0$ is the class of all projective modules in the case when $\mathcal C = \{ R \}$.  

The class $\dFilt{\mathcal C}$ is always closed under arbitrary direct sums and direct summands, but also under (transfinite) extensions. The latter fact is immediate from the following description of $\dFilt{\mathcal C}$ obtained in \cite{ST} (see also \cite[7.12]{GT}): 

\begin{lemma}\label{kapl} Let $R$ be a ring, $\kappa$ be a regular uncountable cardinal, and $\mathcal C$ be a set of $< \kappa$-presented modules. Denote by $\mathcal A$ the class of all $< \kappa$-presented modules from $\dFilt{\mathcal C}$. Then $\dFilt{\mathcal C} = \Filt{\mathcal A}$. In particular, $\dFilt{\mathcal C}$ is $\kappa$-deconstructible.        
\end{lemma}

In the particular case when $\mathcal C = \{ R \}$, Lemma \ref{kapl} says that the class $\dFilt{\mathcal C} = \mathcal P_0$ of all projective modules is $\aleph_1$-deconstructible. This is just the classic Kaplansky Theorem on projective modules \cite[26.2]{AF}: \emph{Each projective module is a direct sum of countably presented projective modules.} 

In view of Lemma \ref{kapl}, our study of deconstructible classes of modules below will include also the classes of the form $\dFilt{\mathcal C}$ for any set of modules $\mathcal C$. 

\medskip
Deconstructible classes of modules are almost omnipresent in homological algebra. For example, for each $n < \omega$, the classes $\mathcal P_n$, $\mathcal F_n$, and $^\perp \mathcal I_n$ are deconstructible, cf.\ \cite[\S 8.1]{GT}. Here, $\mathcal P_n$, $\mathcal F_n$, and $\mathcal I_n$ denote the classes of all modules of projective, flat, and injective, dimension $\leq n$, respectively (so in particular, $\mathcal P_0$, $\mathcal F_0$, and $\mathcal I_0$ are the classes of all projective, flat, and injective modules, respectively). Also $^\perp \mathcal I_n := \{ M \in \rmod R \mid \Ext 1RMN = 0 \hbox{ for each } N \in \mathcal I _n \}$, where Ext$^1$ denotes the first derived functor of Hom; in particular, $\Ext 1RMN = 0$ is equivalent to the splitting of all short exact sequences in $\rmod R$ of the form $0 \to N \to P \to M \to 0$.        

\medskip

For a module $M$, we will denote by $L(M)$ the (complete, modular) lattice of all submodules of $M$, and by $\gen M$ the minimal number of $R$-generators of $M$. A nice feature of the modules $M$ possessing a $\mathcal C$-filtration is that any such filtration can be expanded to a complete distributive sublattice, $\mathcal H$,  of the complete modular lattice $L(M)$ with several remarkable properties. This fact is known as the Hill Lemma, cf.\ \cite[\S7.1]{GT}:

\begin{lemma} \label{Hill} 
Let $R$ be a ring, $\kappa$ an infinite regular cardinal, and $\mathcal C$ a set of $<\kappa$-presented modules. Let $M$ be a module possessing a $\mathcal C$-filtration $\mathcal M = (M_\alpha \mid
\alpha \leq \sigma)$. Then there is a family $\mathcal H$ consisting of submodules of $M$ such that
\begin{itemize}
\item[(H1)] $\mathcal M \subseteq \mathcal H$.
\item[(H2)] $\mathcal H$ is closed under arbitrary sums and intersections. $\mathcal H$ is a complete distributive sublattice of $L(M)$.

\item[(H3)] Let $N,P \in \mathcal H$ be such that $N \subseteq P$. Then the module $P/N$ is $\mathcal C$-filtered.
\item[(H4)] Let $N \in \mathcal H$ and $X$ be a subset of $M$ of cardinality $<\kappa$.  Then there is a $P \in \mathcal H$, such that $N \cup X \subseteq P$ and $P/N$ is $<\kappa$-presented.
\end{itemize}
\end{lemma}

Hill Lemma has a number of strong consequences for the structure theory of modules (see \cite[\S7]{GT} for some of them, including a proof of Lemma \ref{kapl}). Here, we will need a simple corollary that makes it possible to replace the original $\mathcal C$-filtration of a module $M$ by a different one that better fits a particular presentation of $M$:
\begin{corollary}\label{changefilt} Let $R$ be a ring, $\kappa$ an infinite regular cardinal, and $\mathcal C$ a set of $<\kappa$-presented modules. Let $\mathcal C^\prime$ denote the class of all modules possessing a $\mathcal C$-filtration of length $< \kappa$ (whence $\Filt{\mathcal C^\prime} = \Filt{\mathcal C}$).

Let $M$ be a $\mathcal C$-filtered module and $\{ m_\beta \mid \beta < \mu \}$ be a set of $R$-generators of $M$. 

Then there exists a $\mathcal C^\prime$-filtration $\mathcal M^\prime = ( M^\prime_\beta \mid \beta \leq \mu )$ of $M$ such that $m_\beta \in M^\prime_{\beta + 1}$ for each $\beta < \mu$. Moreover, if $\mu = \gen(M) \geq \kappa$ then $\gen(M^\prime_\beta) < \mu$ for each $\beta < \mu$.
\end{corollary}
\begin{proof} Let $\mathcal M = (M_\alpha \mid \alpha \leq \sigma)$ be a $\mathcal C$-filtration of $M$. Consider the corresponding family $\mathcal H$ from Lemma \ref{Hill}. 

The $\mathcal C^\prime$-filtration $\mathcal M^\prime$ will be selected from $\mathcal H$ by induction as follows: $M^\prime_0 = 0$; if $M^\prime_\beta \in \mathcal H$ is defined and $m_\beta \in M^\prime_\beta$, then we put $M^\prime_{\beta + 1} = M^\prime_\beta$. If $m_\beta \notin M^\prime_\beta$, we use property (H4) from Lemma \ref{Hill} to find a module $M^\prime_{\beta +1} \in \mathcal H$ such that $M^\prime_\beta \subseteq M^\prime_{\beta + 1}$, $m_\beta \in M^\prime_{\beta + 1}$, and $M^\prime_{\beta + 1}/M^\prime_\beta$ is $< \kappa$-presented. By property (H3), $M^\prime_{\beta + 1}/M^\prime_\beta$ has a $\mathcal C$-filtration of length $< \kappa$, so $M^\prime_{\beta + 1}/M^\prime_\beta \in \mathcal C^\prime$. 

If $\beta \leq \mu$ is a limit ordinal, we let $M^\prime_\beta = \bigcup_{\gamma < \beta} M^\prime_\gamma$, which is a module from $\mathcal H$ by property (H2). Since $m_\beta \in M^\prime_{\beta + 1}$ for each $\beta < \mu$, $M^\prime_\mu = M$.  

Finally, for each $\beta < \mu$, $(M^\prime_\gamma \mid \gamma \leq \beta)$ is a $\mathcal C^\prime$-filtration of $M^\prime_\beta$ of length $\beta$. So if $\mu = \gen(M) \geq \kappa$, then $\gen(M^\prime_\beta) < \mu$ because $\kappa$ is regular.    
\end{proof}

\subsection{Categoricity} 
We will consider the following notions of categoricity for a class of modules $\mathcal D$ (cf.\ \cite{M2}): If $\lambda$ is an infinite  cardinal then $\mathcal D$ is \emph{$\lambda$-categorical} (or \emph{categorical in $\lambda$}) provided that $\mathcal D$ contains a module of cardinality $\lambda$, and all modules of cardinality $\lambda$ contained in $\mathcal D$ are isomorphic. The class $\mathcal D$ is \emph{categorical in a tail} of cardinals, provided there exists an infinite cardinal $\kappa$ such that $\mathcal D$ is \emph{$\lambda$-categorical} for each cardinal $\lambda \geq \kappa$. 

Let $\mathcal D$ be a deconstructible class of modules, and let $\kappa$ be the least cardinal $\kappa \geq \card R + \aleph_0$ such that $\mathcal D$ is $\kappa^+$-deconstructible. The \emph{categoricity spectrum} of $\mathcal D$ is defined as the class of all cardinals $\lambda \geq \kappa$ such that $\mathcal D$ is $\lambda$-categorical. This spectrum is denoted by $\Spec {\mathcal D}$. We will say that $\Spec {\mathcal D}$ is \emph{bounded} in case there exists a cardinal $\mu \geq \kappa$ such that $\Spec {\mathcal D} \subseteq \langle \kappa , \mu \rangle$. For an infinite cardinal $\nu$, we will denote by $\langle \nu, \infty )$ the class of all cardinals $\geq \nu$.   

\medskip
Obviously, categoricity in a tail implies categoricity in a big cardinal, but there are various instances where also the reverse implication holds, see \cite{M1}, \cite{M2} for some recent results related to this phenomenon. In classic model theory of first-order theories, $T$, the phenomenon is expressed by the Morley-Shelah Categoricity Theorem: if $T$ is categorical in a cardinal $\kappa > \mid{T}\mid$, then $T$ is categorical in all cardinals $\lambda > \mid{T}\mid$, \cite{S}. Here $\mid{T}\mid$ denotes the cardinality of non-logical symbols that occur in the axioms of $T$. Notice however, that classes of modules of interest in the present paper, i.e., those of the form $\Filt{\mathcal C}$ and $\dFilt{\mathcal C}$ for a set of modules $\mathcal C$, are not first-order axiomatizable in general. For example, the class $\mathcal F _0$ of all flat modules is first-order axiomatizable, iff the ring $R$ is left coherent, while the class $\mathcal P _0$ of all projective modules is first-order axiomatizable, iff $R$ is left coherent and right perfect, cf.\ \cite{ES}.  

The Morley-Shelah Categoricity Theorem is one of the major achievements of classic model theory. So it is natural to ask for its extension to the setting of infinitary model theory. Here, we will focus on the main setting of contemporary infinitary model theory, the abstract elementary classes in the sense of \cite{S-aec1} and \cite{S-aec2} (see also \cite{B}). First, we recall the relevant definitions:

\begin{definition} \label{aec}
A pair $(\mathcal A, \preceq )$ is an \emph{abstract elementary class} (or an \emph{AEC}), if $\mathcal A$ is a class of structures (in a fixed vocabulary $\tau$), and $\preceq$ is a partial order on $\mathcal A$, both $\mathcal A$ and $\preceq$ are closed under isomorphisms, and satisfy the following axioms (A1)-(A4):
\begin{itemize}
\item[{\rm{(A1)}}] If $A \preceq B$, then $A$ is a substructure of $B$.
\item[{\rm{(A2)}}] If $( A_i \mid i < \delta )$ is a continuous $\preceq$-increasing chain of elements of $\mathcal A$, then
\begin{enumerate}
\item $\bigcup_{i<\delta}A_i \in \mathcal A$;
\item $A_j \preceq \bigcup_{i<\delta}A_i$ for each $j<\delta$;
\item If $M \in \mathcal A$ and $A_i \preceq M$ for each $i<\delta$, then $\bigcup_{i<\delta}A_i \preceq M$.
\end{enumerate}
\item[{\rm{(A3)}}] If $A,B,C \in \mathcal A$, $A \preceq C$, $B \preceq C$ and $A$ is a substructure of $B$, then $A \preceq B$.
\item[{\rm{(A4)}}] There is a cardinal $\kappa$ such that, if $A$ is a subset of $B \in \mathcal A$, then there is
$A^\prime \in \mathcal A$ which contains $A$, such that $A^\prime \preceq B$ and the cardinality of
$A^\prime$ is at most $\card{A} + \kappa$. The least such infinite cardinal $\kappa \geq \card \tau$ is called the \emph{L\"owenheim-Skolem number}
of $(\mathcal A, \preceq )$ and denoted by $\hbox{LS}(\mathcal A)$.
\end{itemize}
Here, $\preceq$ is said to be \emph{closed under isomorphisms} provided that $N^\prime \preceq M^\prime$, whenever $N \preceq M$ and there is an isomorphism
$f : M \to M^\prime$, such that $f(N) = N^\prime$.

Also, $( A_i \mid i< \delta )$ is a \emph{continuous $\preceq$-increasing chain} provided that
\begin{itemize}
\item[{\rm{(i)}}]  $A_i \in \mathcal A$,
\item[{\rm{(ii)}}] $A_i \preceq A_{i+1}$ for all $i < \delta$, and
\item[{\rm{(iii)}}] $A_i = \bigcup_{j<i}A_j$ for all limit ordinals $i < \delta$.
\end{itemize}
\end{definition}

Basic examples of AECs are $(\mathcal A _T, \preceq _T)$, where $\mathcal A _T$ is the class of all models of a first-order theory $T$, and $\preceq _T$ is the relation of being a submodel, or the relation of being an elementary submodel. In the particular setting of the first-order theory of modules over a ring $R$, the term submodel just means (right $R$-) submodule. In this setting, also the pair $(\rmod R, \preceq^*)$ is an AEC, where $\preceq^*$ is the relation of being a pure submodule.  

\medskip
Deconstructible classes of modules are sources of particular kinds of AEC's of modules, called the AEC's of roots of Ext. The following definition and theorem go back to \cite{BET} (see also \cite[\S 10.3]{GT}): 

\begin{definition}\label{AECroots} Let $R$ be a ring, $\mathcal C$ a class of modules, and let 
$$\mathcal A = {}^{\perp_\infty} \mathcal C := \{ M \in \rmod R \mid \Ext iRMC = 0 \hbox{ for all } i \geq 1 \hbox{ and } C \in \mathcal C \}.$$ 
Consider the pair $(\mathcal A, \preceq )$, where $A \preceq B$, iff $A$ is a submodule of $B$ such that $A, B, B/A \in \mathcal A$.
\end{definition}

\begin{theorem}\label{bet-cluj} Let $R$ be a ring, $\mathcal C$ a class of modules, $\mathcal A = {}^{\perp_\infty} \mathcal C$, and $(\mathcal A, \preceq )$ be as in Definition \ref{AECroots}. Let $\kappa$ be an infinite cardinal $\geq \card R$. Then the following conditions are equivalent:
\begin{itemize}
\item[{\rm{(a)}}] $(\mathcal A, \preceq )$ is an AEC with $\hbox{LS}(\mathcal A) \leq \kappa$.
\item[{\rm{(b)}}] $\mathcal A$ is a $\kappa^+$-deconstructible class closed under direct limits.
\end{itemize}
In particular, $(\mathcal A, \preceq )$ is an AEC, if and only if $\mathcal A$ is a deconstructible class closed under direct limits. \footnote{\it Added in proof: \rm Recently, Jan \v Saroch and the author proved that if $(\mathcal A, \preceq )$ is an arbitrary AEC such that $\mathcal A$ is a deconstructible class of modules and the relation $\preceq$ refines the direct summand relation, then $\mathcal A$ is closed under direct limits, see \cite{SaT}.}
\end{theorem}

The AEC's from Theorem \ref{bet-cluj}, i.e., those of the form $(\mathcal A, \preceq)$, where $\mathcal A = {}^{\perp_\infty} \mathcal C$ for a class $\mathcal C \subseteq \rmod R$,  are called the \emph{AEC's of roots of Ext}. As the underlying structures for $\mathcal A$ are (right $R$-) modules, we always have $\hbox{LS}(\mathcal A) \geq \card R + \aleph_0$ by Definition \ref{aec}.(A4). 

\medskip
A long-standing open problem of abstract model theory concerns categoricity for the general abstract elementary classes defined in \ref{aec}, see e.g.\ \cite[vol.1, N.4.3]{S-aec2}. Assume that $(\mathcal A, \preceq)$ is an AEC with L\"{o}wenheim-Skolem number $\kappa$. 

\medskip
{\bf Shelah's Categoricity Conjecture} (SCC) says that if $\mathcal A$ is categorical in a cardinal $\lambda \geq \beth_{(2^\kappa)^+}$, then $\mathcal A$ is categorical in all cardinals $\geq \beth_{(2^\kappa)^+}$. 

\medskip
Here, the function $\beth$ (beth) is defined by induction on the set of all ordinals as follows: $\beth_0 = \aleph_0$, $\beth_{\alpha + 1} = 2^{\beth_{\alpha}}$, and $\beth_\alpha = \sup_{\beta < \alpha} \beth_{\beta}$ when $\alpha$ is a limit ordinal. So for each ordinal $\alpha$, $\aleph_{\alpha} \leq \beth_\alpha$, and the two functions coincide, iff GCH holds.

\medskip
The two main results of the present paper are as follows: in Theorem \ref{roots} of Section \ref{projaec}, we prove SCC for all AEC's of roots of Ext. In Section \ref{cat}, we turn to general deconstructible classes of modules. In Theorem \ref{main}, we show that the categoricity of $\mathcal D$ in a big cardinal implies its categoricity in a tail of cardinals, for an arbitrary deconstructible class of modules $\mathcal D$.    

\section{Categoricity for projective modules and the AEC's of roots of Ext}\label{projaec}

Recall that for any ring $R$, the \emph{rank} of a free module $F$ is a cardinal $\kappa$ such that $F \cong R^{(\kappa)}$. Clearly, if $F$ is not finitely generated, then the rank of $F$ is uniquely determined by $F$. It follows that the class of all free modules is $\lambda$-categorical for each cardinal $\lambda > \card R + \aleph_0$. 

For projective modules, we have the following characterization, answering partially Question 4.2 from \cite{M2}:      

  \begin{proposition}\label{projSCC}
	Let $R$ be a ring and $S$ be a representative set of all countably generated projective modules. Let $\kappa = \card R + \aleph_0$ and $\nu = \kappa + \card S$. Then $\nu \leq 2^\kappa$, and the following conditions are equivalent:
	\begin{enumerate}
		\item[(1)] The class $\mathcal P _0$ is $\lambda$-categorical for all cardinals $\lambda > \nu$.
		\item[(2)] $\mathcal P _0$ is $\lambda$-categorical for some cardinal $\lambda \geq \kappa$.
		\item[(3)] The module $P^{(\omega)}$ is free, for each countably generated projective module $P$.
	\end{enumerate}
   \end{proposition}
\begin{proof} Since $\card {R^{(\omega)}} = \kappa$, we have $\card S \leq 2^\kappa$, and hence also $\nu \leq 2^\kappa$. Clearly (1) implies (2). 

Assume (2) and let $P$ be a non-zero countably generated projective module. Since $\lambda \geq \kappa$, both $R^{(\lambda)}$ and $P^{(\lambda)}$ have cardinality $\lambda$, whence (2) yields an isomorphism $P^{(\lambda)} \cong R^{(\lambda)}$. Consequently, the countably generated module $P^{(\omega)}$ is isomorphic to a direct summand in $R^{(\omega)}$, i.e., $P^{(\omega)} \oplus X \cong R^{(\omega)}$ for a module $X$. Similarly, $R^{(\omega)} \oplus Y \cong P^{(\omega)}$ for a module Y. Then
$$P^{(\omega)} \cong (R^{(\omega)} \oplus R^{(\omega)}) \oplus Y \cong R^{(\omega)} \oplus P^{(\omega)}$$
and similarly $R^{(\omega)} \oplus P^{(\omega)} \cong R^{(\omega)}$. This proves (3).

Assume (3). Let $P$ be a projective module of cardinality $\lambda > \nu$. By the classic Kaplansky Theorem \cite[26.2]{AF}, $P$ is isomorphic to a direct sum of copies of modules from the set $S$, i.e., $P \cong \bigoplus_{Q \in S} Q^{(\kappa_Q)}$ for some cardinals $\kappa_Q \geq 0$ ($Q \in S$). Let $P^\prime = \bigoplus_{Q \in S^\prime} Q^{(\kappa_Q)}$, where $S^\prime$ is the subset of $S$ consisting of all the modules $Q$ with  $\kappa_Q < \aleph_0$. By (3), for each $Q \in S \setminus S^\prime$, we have $Q^{(\kappa_Q)} \cong R^{(\kappa_Q)}$. Since $\card P = \lambda > \card S \geq \card S^\prime$, we infer that $\card {P^\prime} < \lambda$, and $P \cong P^\prime \oplus \bigoplus_{Q \in S \setminus S^\prime} R^{(\kappa_Q)}$. Thus the free module $\bigoplus_{Q \in S \setminus S^\prime} R^{(\kappa_Q)}$ has cardinality $\lambda$, so $P \cong P^\prime \oplus R^{(\lambda)}$. By Eilenberg's Trick (see e.g.\ \cite[Proposition 18.1]{P}), $P$ is isomorphic to $R^{(\lambda)}$. Thus, the class $\mathcal P _0$ is categorical in $\lambda$, and (1) holds true.
\end{proof}

\begin{remark}\label{Bass} If $R$ is a ring such that all projective modules are free, then Condition (3) above obviously holds true. By a classic result of Bass \cite[Theorem 18.7]{P}, countably generated projective modules over any commutative noetherian domain $R$ are free, so Condition (3) also holds in this case. However, (3) fails for any ring $R$ containing a non-trivial central idempotent. 
\end{remark}  

Next, we show that SCC holds for all AEC's of roots of Ext:

\begin{theorem}\label{roots} Let $R$ be a ring and $(\mathcal D, \preceq )$ be an AEC of roots of Ext. Then SCC holds for $(\mathcal D, \preceq )$.  
\end{theorem}
\begin{proof} Let $\rho = \hbox{LS}(\mathcal D)$. By Theorem \ref{bet-cluj}, $\mathcal D$ is $\rho^+$-deconstructible. Moreover, $\mathcal D$ contains all projective modules, and $\mathcal D$ is closed under direct limits. 

First, assume that $\mathcal D$ contains a non-projective module. Since $\mathcal D$ is $\rho^+$-deconstructible, $\mathcal D$ contains a $\leq \rho$-presented non-projective module $N$. Then for each cardinal $\lambda \geq \rho$, $N^{(\lambda)}$ and $R^{(\lambda)}$ are non-isomorphic modules, both of cardinality $\lambda$. This proves that $\mathcal D$ is not $\lambda$-categorical for any $\lambda \geq \rho$. Since $\rho \leq \beth_\rho < \beth_{(2^\rho)^+}$, SCC trivially holds for $(\mathcal D, \preceq )$. 

In the remaining case of $\mathcal D = \mathcal P _0$, the class of all projective modules is closed under direct limits, so $R$ is a right perfect ring. Hence $R$ contains a finite basic set of primitive idempotents $\{ e_1, \dots , e_n \}$, and up to isomorphism, each projective module is of the form $\bigoplus_{0 < i \leq n} (e_iR)^{(\kappa_i)}$ for a unique $n$-tuple of cardinals $(\kappa_1, \dots , \kappa_{n})$, \cite[27.11]{AF}. In particular, the set $S$ from Proposition \ref{projSCC} is countable, so (in the notation of Proposition \ref{projSCC}) $\nu = \kappa$ and $\kappa = \rho$, so $\kappa^+ = \rho^+ < \beth_{(2^\rho)^+}$, and SCC follows by Proposition \ref{projSCC}. 
\end{proof}

\begin{remark}\label{stronger} The proof of Theorem \ref{roots} actually shows that in the case of the AEC's of roots of Ext, SCC holds in a stronger form, with $\beth_{(2^\rho)^+}$ replaced by $\rho ^+ = \hbox{LS}(\mathcal D)^+$. Indeed, in the first case of the proof, $\Spec {\mathcal D} = \emptyset$, while in the second case either $\Spec {\mathcal D} \supseteq \langle \rho{^+}, \infty )$ when $n = 1$ (i.e., when $R$ is isomorphic to a full matrix ring over a local right perfect ring), or else $\Spec {\mathcal D} = \emptyset$ (when $n > 1$). 
\end{remark}  

\section{Categoricity for deconstructible classes of modules}\label{cat}

In order to proceed to arbitrary deconstructible classes of modules, we recall the notion of a strong splitter that generalizes the notion of a projective module: 

\begin{definition}\label{ssplit}
A module $M$ is a \emph{splitter} (or an \emph{exceptional module}) provided that $\Ext 1RMM = 0$. Moreover, $M$ is a \emph{strong splitter}, if $\Ext 1RM{M^{(I)}} = 0$ for any set $I$. 
\end{definition}

Notice that if $M$ is a strong splitter, then so is any direct summand of $M$, as well as $M^{(\kappa)}$ for any cardinal $\kappa$.  

Strong splitters are quite common in module theory. For example, each (infinitely generated) tilting module is a strong splitter, as the defining condition of a strong splitter from Definition \ref{ssplit} is just the condition (T2) from the definition of a tilting module (see e.g.\ \cite[13.1]{GT}). Also notice that it suffices to test this condition only at a (suitable) single set $I$:

\begin{lemma}\label{single}
Let $R$ be a ring, $\kappa \geq \aleph_0$ be a cardinal, and $M$ be a $\leq \kappa$-presented module. Then $M$ is a strong splitter, iff $\Ext 1RM{M^{(\kappa)}} = 0$.  

Moreover, every finitely presented splitter is strong.
\end{lemma}
\begin{proof} By assumption, there is a presentation $0 \to K \to F \to M \to 0$ of $M$ such that $F$ and $K$ are $\leq \kappa$ generated modules, and $F$ is free. 

Let $I$ be a set of cardinality $> \kappa$ and $\varphi \in \Hom RK{M^{(I)}}$. Then there is a subset $J \subseteq I$ of cardinality $\leq \kappa$ such that $\im \varphi \subseteq M^{(J)}$. 
If $\Ext 1RM{M^{(\kappa)}} = 0$, then $\varphi$ extends to some $\psi \in \Hom RF{M^{(J)}} \subseteq \Hom RF{M^{(I)}}$. Thus also $\Ext 1RM{M^{(I)}} = 0$.     

For the moreover part, we proceed similarly, taking a presentation of $M$ with $K$ and $F$ finitely generated, using the facts that $\im \varphi \subseteq M^{(J)}$ for a finite subset $J$ of $I$ and $\Ext 1RM{M^{(J)}} = 0$.  
\end{proof}

Surprisingly, each deconstructible class of modules contains a proper class of splitters: 

\begin{lemma}\label{et} Let $R$ be a ring, $\mathcal D$ be a deconstructible class of modules, and $0 \neq D \in \mathcal D$. Let 
$\kappa = \card D + \card R + \aleph_0$, and let $\lambda > \kappa$ be such that $\lambda ^\kappa = \lambda$  (e.g., $\lambda = 2^\rho$ for some $\rho \geq \kappa$). Then there exists a $\{ D \}$-filtered splitter $N_\lambda$ of cardinality $\lambda$, such that $D \subseteq N_\lambda$, $\Ext 1RD{N_\lambda} = 0$, and $N_\lambda \in \mathcal D$.
\end{lemma}
\begin{proof} By \cite[Theorem 2]{ET}, there exists a module $N_\lambda$ of cardinality $\lambda$ such that $D \subseteq N_\lambda$, $\Ext 1RD{N_\lambda} = 0$, and $N_{\lambda}/D$ is $\{ D \}$-filtered. Then also $N_\lambda$ is $\{ D \}$-filtered, whence $N_\lambda \in \mathcal D$. Since $\Ext 1RD{N_\lambda} = 0$, also $\Ext 1R{N_\lambda}{N_\lambda} = 0$ by 
the Eklof Lemma \cite[6.2]{GT}.  
\end{proof} 

By Lemma \ref{single}, all finitely presented splitters in a deconstructible class $\mathcal D$ are strong splitters. However, despite Lemma \ref{et}, $\mathcal D$ need not contain any non-zero strong splitters in general:

\begin{example}\label{vonN} Let $R$ be a simple right hereditary ring containing an infinite set of non-zero pairwise orthogonal idempotents $\{ e_i \mid i < \omega \}$. For a concrete example, consider $R = \varinjlim_{n < \omega} M_{2^n}(K)$, the direct limit of the $\omega$-directed system of full $2^n \times 2^n$-matrix rings $M_{2^n}(K)$ over a fixed field $K$, and morphisms induced by the block diagonal ring embeddings $M_{2^n}(K) \to M_{2^{n+1}}(K)$. 

Let $N = R/(\bigoplus_{i < \omega} e_iR)$. Then $\Ext 1RN{P^{(\omega)}} \neq 0$ for any non-zero module $P$, because $Re_iR = R$, whence $Pe_i \neq 0$, for each $i < \omega$, so there exist homomorphisms $f \in \Hom R{\bigoplus_{i < \omega} e_iR}{P^{(\omega)}}$ that cannot be extended to $R$.

Let $M$ be any module containing $N$, and $\mathcal D = \Filt {\{ M \}}$. Since $R$ is right hereditary, \emph{each} non-zero module $Q \in \mathcal D$ satisfies $\Ext 1RQ{P^{(\omega)}} \neq 0$ for any non-zero module $P$, because $Q$ contains a copy of $N$. In particular, $Q$ is not a strong splitter.
\end{example}

\medskip
By Proposition \ref{projSCC}, categoricity of the class $\mathcal D = \mathcal P _0$ in a tail is equivalent to its categoricity in a big cardinal. The same holds for any deconstructible class $\mathcal D$ which contains a non-zero strong splitter: 

  \begin{lemma}\label{splitSCC}
Let $R$ be a ring, and $\mathcal D$ be a deconstructible class of modules containing a non-zero strong splitter $M$. Let $\kappa$ be an infinite cardinal $\geq \card R + \card M$ such that $\mathcal D$ is $\kappa^+$-deconstructible. Then for each $\lambda > \kappa$, the following conditions are equivalent:
	\begin{enumerate}
		\item[(1)] The class $\mathcal D$ is $\mu$-categorical for each cardinal $\mu \geq \lambda$.
		\item[(2)] $\mathcal D$ is $\lambda$-categorical.
		\item[(3)] For each cardinal $\mu \geq \lambda$, all the elements of $\mathcal D$ of cardinality $\mu$ are isomorphic to $M^{(\mu)}$.
	\end{enumerate}
   \end{lemma}
\begin{proof} We only have to prove that (2) implies (3). Assume (2) and let $\mu$ be a cardinal $\geq \lambda$. Let $N \in \mathcal D$ be of cardinality $\mu$. If $\mu = \lambda$, then $N \cong M^{(\mu)}$ by (2). 

Assume $\mu > \lambda$ and let $\mathcal N = ( N_\beta \mid \beta \leq \sigma )$ be a $\mathcal C$-filtration of $N$. Since $\gen N = \mu$, and all the modules in $\mathcal C$ have cardinality $\leq \kappa < \mu$, necessarily $\sigma \geq \mu$. By Corollary \ref{changefilt}, possibly replacing $\mathcal C$ with the class $\mathcal C^\prime$ of all modules possessing a $\mathcal C$-filtration of length $\leq \kappa$, we can w.l.o.g.\ assume that $\sigma = \mu$.

Let $\alpha$ be the minimal element of the set of all ordinals $\beta \leq \mu$ such that $N_\beta$ has cardinality $\geq \lambda$. Since $\mathcal C$ consists of $\leq \kappa$-generated modules, and $\kappa < \lambda$, $\alpha$ is a limit ordinal. By the continuity of $\mathcal N$, $\card {N_\alpha} = \sup_{\gamma < \alpha} \card {N_\gamma}$. By the minimality of $\alpha$, the latter supremum is $\leq \lambda$, so $\card {N_\alpha} = \lambda$. By (2), $N_\alpha \cong M^{(\lambda)}$.

Since $\alpha < \mu$, we can proceed by transfinite induction, replacing $N$ by $N/N_\alpha$ and $\mathcal N$ by $\bar{\mathcal N} = ( N_\beta/N_\alpha \mid \alpha \leq \beta \leq \sigma )$ etc., and select from $\mathcal N$ a $\mathcal E$-filtration $( E_\gamma \mid \gamma \leq \mu )$ of $N$ where $\mathcal E = \{ M^{(\lambda)} \}$. Since $M$ is a strong splitter, by induction on $\gamma < \mu$, we infer that $E_{\gamma} \cong M^{(\lambda . \gamma)}$ and $E_{\gamma +1} \cong E_{\gamma} \oplus M^{(\lambda)}$ for each $\gamma < \mu$. Thus $N \cong M^{(\mu)}$, and (3) holds.         
\end{proof}

\begin{corollary}\label{tailcor}
Let $R$ be a ring, and $\mathcal D$ be a deconstructible class of modules containing a non-zero strong splitter $M$. Let $\mu$ be the least infinite cardinal $\geq \card R$ such that $\mathcal D$ is $\mu^+$-deconstructible. Let $\kappa$ be the least cardinal such that $\mu + \card M \leq \kappa$. 

Then either $\Spec {\mathcal D} \subseteq \langle \mu , \kappa \rangle$, or there is a cardinal $\lambda > \kappa$ such that $\Spec {\mathcal D} \cap ( \kappa , \infty ) = \langle \lambda , \infty )$. 
\end{corollary}
\begin{proof} This follows by Lemma \ref{splitSCC}: if $\Spec {\mathcal D} \nsubseteq \langle \mu , \kappa \rangle$, we take the least $\lambda > \kappa$ such that $\mathcal D$ is $\lambda$-categorical. 
\end{proof}

We pause to consider an example where categoricity in a tail is easily related to the structure of the ring $R$:

\begin{example}\label{pinj} Let $R$ be a ring and $\kappa = \card R + \aleph_0$. Let $\mathcal I$ be any class of pure-injective modules and $\mathcal D = {}^{\perp_{\infty}} \mathcal I = \{ M \in \rmod R \mid \Ext iRMI = 0 \hbox{ for all } I \in \mathcal I \hbox{ and } i \geq 1 \}$. Denote by $\mathcal I ^\prime$ the class of all cosyzygies of modules from $\mathcal I$ (cf.\ \cite[Corollary 2.25]{GT}). Then $\mathcal D = {}^\perp \mathcal I ^\prime$, and $\mathcal I ^\prime$ consists of pure-injective modules by \cite[Lemma 6.20]{GT}. By \cite[Lemma 6.17]{GT}, the class $\mathcal D$ is $\kappa^+$-deconstructible for $\kappa = \card R + \aleph_0$, and $\mathcal D$ is closed under direct limits by \cite[Corollary 2.9]{GT}. Hence $(\mathcal D, \preceq )$ is an AEC of roots of Ext by Theorem \ref{bet-cluj}. Since $R \in \mathcal D$, Lemma \ref{splitSCC} applies for $M = R$. Notice that $\mathcal D$ contains the class $\mathcal F _0$ of all flat modules.

For $\lambda > \kappa$, Condition (3) of Lemma \ref{splitSCC} says that all modules in $\mathcal D$ of cardinality $\geq \lambda$ are free. Hence all flat modules are projective, and $R$ is a right perfect ring. Moreover, by \cite[27.11]{AF}, $R = (eR)^n$ for a primitive idempotent $e \in R$, that is, $R \cong M_n(eRe)$ is a full matrix ring over the local right perfect ring $eRe$. 

If $\mathcal D = \mathcal F_0$, then it is easy to see that also the converse holds true, as the $\lambda$-categoricity of $\mathcal F_0$ holds for each $\lambda > \kappa$ when $R$ is a full matrix ring over a local right perfect ring. Thus we recover the characterization given in \cite[Theorem 3.26]{M2}. 

If $\mathcal D \supsetneq \mathcal F_0$, then $\mathcal D$ contains a $\leq \kappa$-generated non-flat module, whence -- for any ring $R$ -- the $\lambda$-categoricity of $\mathcal D$ fails for all $\lambda \geq \kappa$.            
\end{example}

\begin{remark}\label{add} Condition (3) from Lemma \ref{splitSCC} implies that all big modules in the class $\mathcal D$ are isomorphic to direct sums of copies of the module $M$. In fact, if $\Add M$ denotes the class of all direct summands of arbitrary direct sums of copies of the module $M$, then necessarily $\mathcal D \subseteq \Add M$. 

Indeed, if $N \in \mathcal D$, then Condition (3) implies that a sufficiently big direct sum of copies of the module $N$ is isomorphic to a direct sum of copies of the module $M$, whence $N \in \Add M$. 
\end{remark}

\medskip
The picture changes completely when $\mathcal D$ contains a module which is not a strong splitter:

\begin{lemma}\label{notail} Let $R$ be a ring, and $\mathcal D$ be a deconstructible class of modules containing a module $M$ which is not a strong splitter. Let $\kappa$ be an infinite cardinal $\geq \card R + \card M$. 

Assume that $\mathcal D$ is $\lambda$-categorical for some $\lambda \geq \kappa$. Then $\mathcal D$ is not $\mu$-categorical for \emph{any} $\kappa \leq \mu$ such that $\mu \neq \lambda$.   

Hence $\card \hbox{Spec}_{\kappa}(\mathcal D) \leq 1$, where $\hbox{Spec}_{\kappa}(\mathcal D) = \Spec {\mathcal D} \cap \langle \kappa , \infty )$. 
\end{lemma}  
\begin{proof} Clearly, it suffices to prove the claim under the additional assumption that $\lambda$ is the least cardinal $\geq \kappa$ such that $\mathcal D$ is $\lambda$-categorical. Then $\mathcal D$ is not $\mu$-categorical for each $\kappa \leq \mu < \lambda$, and 
all modules of cardinality $\lambda$ in the class $\mathcal D$ are isomorphic to $M^{(\lambda)}$. 

By Lemma \ref{single}, $\Ext 1RM{M^{(\lambda)}} \neq 0$, so there is a non-split short exact sequence of the form 
$$ (*) \quad 0 \to M^{(\lambda)} \to M^{(\lambda)} \to M \to 0.$$ 
By induction, define a strictly increasing chain of modules $( M_\alpha \mid \alpha < \lambda^+ )$ such that $M_\alpha \cong M^{(\lambda)}$ for all $\alpha < \lambda ^+$ as follows: $M_0 = M^{(\lambda)}$. If $M_\alpha$ is defined, then we use the short exact sequence (*) to extend $M_\alpha$ to a module $M_{\alpha +1} \in \mathcal D$ such that $M_{\alpha +1}/M_\alpha \cong M$, and the embedding $M_\alpha \subseteq M_{\alpha +1}$ does not split. Then $\card M_{\alpha +1} = \card M_{\alpha} = \card M^{(\lambda)} = \lambda$, so $M_{\alpha +1} \cong M^{(\lambda)}$ because $\mathcal D$ is $\lambda$-categorical. If $\alpha < \lambda^+$ is a limit ordinal, we let $M_\alpha = \bigcup_{\beta < \alpha} M_\beta$. Then $M_\alpha \in \mathcal D$ and $\card M_\alpha = \lambda$, whence again $M_\alpha \cong M^{(\lambda)}$.

Let $N = \bigcup_{\alpha < \lambda^+} M_\alpha$. Then $N \in \mathcal D$. First, we will prove that $N \ncong M^{(\lambda^+)}$; this will imply that the class $\mathcal D$ is not $\lambda^+$-categorical.               

Assume $N \cong M^{(\lambda^+)}$, so there exist submodules $P_\alpha$ ($\alpha < \lambda^+$) of $N$ such that $N = \bigoplus_{\alpha < \lambda^+} P_\alpha$ and $P_\alpha \cong M$ for each $\alpha < \lambda^+$. Let $\alpha_0 = 0$ and $I_0$ be the least subset of $\lambda^+$ such that $M_{\alpha_0} \subseteq \bigoplus_{\beta \in I_0} P_\beta$. Then $I_0$ has cardinality $\leq \lambda$ and there exists $\alpha_0 \leq \alpha_1 < \lambda^+$ such that $\bigoplus_{\beta \in I_0} P_\beta \subseteq M_{\alpha_1}$. Proceeding in this way, we obtain a chain of ordinals $0 = \alpha_0 \leq \alpha_1 \leq \dots \leq \alpha_n \leq \alpha_{n+1} \leq \dots$, and a chain of subsets $I_0 \subseteq I_1 \subseteq \dots \subseteq I_n \subseteq I_{n+1} \subseteq \dots$ of $\lambda^+$ of cardinality $\leq \lambda$, such that for each $n < \omega$, $M_{\alpha_n} \subseteq \bigoplus_{\beta \in I_{n}} P_\beta \subseteq M_{\alpha_{n+1}}$. Let $\alpha = \sup_{n < \omega} \alpha_n (< \lambda^+)$ and $I = \bigcup_{n < \omega} I_n$. Then $M_\alpha = \bigcup_{n < \omega} M_{\alpha_n} = \bigoplus_{\alpha \in I} P_\alpha$ is a direct summand in $N$, and hence also in $M_{\alpha + 1}$. This contradicts our definition of the module $M_{\alpha + 1}$.      

The remaining case of $\mu > \lambda^+$ is proved similarly, but the proof is a little more technical: Assume that $N^{(\mu)} \cong M^{(\mu)}$, so there exist submodules $P_\beta$ ($\beta < \mu$) of $N^{(\mu)}$ such that $N^{(\mu)} = \bigoplus_{\beta < \mu} P_\beta$ and $P_\beta \cong M$ for each $\beta < \mu$, and also submodules $N_\gamma$ ($\gamma < \mu$) of $N^{(\mu)}$ such that $N^{(\mu)} = \bigoplus_{\gamma < \mu} N_\gamma$ and $N_\gamma = N$ for each $\gamma < \mu$. For each $\gamma < \mu$, let $\nu_{\gamma} : N \to N^{(\mu)}$ be the $\gamma$th-canonical embedding. 

Let $\alpha_0 = 0$, $J_0 = \{ \alpha_0 \}$, and $I_0$ be the least subset of $\mu$ such that $\nu_0(M_{\alpha_0}) \subseteq \bigoplus_{\beta \in I_0} P_\beta$. Then $I_0$ has cardinality $\leq \lambda$, so there exist $\alpha_0 \leq \alpha_1 < \lambda^+$ and a subset $J_1$ of $\mu$ of cardinality $\leq \lambda$, such that $J_0 \subseteq J_1$ and $\bigoplus_{\beta \in I_0} P_\beta \subseteq \bigoplus_{\gamma \in J_1} \nu_{\gamma}(M_{\alpha_1})$. 

Proceeding in this way, we obtain a chain of ordinals $0 = \alpha_0 \leq \alpha_1 \leq \dots \leq \alpha_n \leq \alpha_{n+1} \leq \dots < \lambda^+$, 
a chain of subsets $J_0 \subseteq J_1 \subseteq \dots \subseteq J_n \subseteq J_{n+1} \subseteq \dots$ of $\mu$ of cardinality $\leq \lambda$, and a chain of subsets $I_0 \subseteq I_1 \subseteq \dots \subseteq I_n \subseteq I_{n+1} \subseteq \dots$ of $\mu$ of cardinality $\leq \lambda$, such that for each $n < \omega$, 
$$\bigoplus_{\gamma \in J_n} \nu_{\gamma}(M_{\alpha_n}) \subseteq \bigoplus_{\beta \in I_{n}} P_\beta \subseteq \bigoplus_{\gamma \in J_{n+1}} \nu_{\gamma}(M_{\alpha_{n+1}}).$$ 
Let $\alpha = \sup_{n < \omega} \alpha_n (< \lambda^+)$, $J = \bigcup_{n < \omega} J_n$, and $I = \bigcup_{n < \omega} I_n$. Then 
$$M^\prime_{\alpha} := \bigoplus_{\gamma \in J} \nu_{\gamma}(M_{\alpha}) = \bigcup_{n < \omega} \bigoplus_{\gamma \in J_n} \nu_{\gamma}(M_{\alpha_n}) = \bigoplus_{\beta \in I} P_\beta$$ is a direct summand in $N^{(\mu)}$, and hence also in $M^\prime_{\alpha + 1} := \bigoplus_{\gamma \in J} \nu_{\gamma}(M_{\alpha + 1})$. This contradicts the fact for each $\gamma \in J$, $\nu_\gamma(M_\alpha)$ is not a direct summand in $\nu_\gamma(M_{\alpha + 1})$ (as $M_{\alpha}$ is not a direct summand in $M_{\alpha + 1}$). This proves that $N^{(\mu)} \ncong M^{(\mu)}$.        

The final claim is just a restatement of the above in terms of properties of $\Spec {\mathcal D}$. 
\end{proof}	

Let $\mathcal D$ be an arbitrary deconstructible class of modules. Consider a subset $\mathcal C \subseteq \mathcal D$ such that $\mathcal D = \Filt{\mathcal C}$. Let $M_{\mathcal C} := \bigoplus_{C \in \mathcal C} C (\in \mathcal D)$. Notice that if $\mathcal D$ is $\mu^+$-deconstructible for an infinite cardinal $\mu \geq \card R$, then we can assume that $\mathcal C$ consists of $\leq \mu$-presented modules, and there are at most $2^\mu$ non-isomorphic modules in $\mathcal C$. Hence we can w.l.o.g.\ assume that $\card {M_{\mathcal C}} \leq 2^\mu$.  

Then there are the following alternatives:

\medskip
{\bf Alternative 1: $M_{\mathcal C}$ is a strong splitter.} In this case, \emph{all} modules in the class $\mathcal D = \Filt{\mathcal C}$ are strong splitters, and they are isomorphic to direct sums of modules from $\mathcal C$. Indeed, let $M \in \mathcal D$, so $M$ has a $\mathcal C$-filtration $( M_\alpha \mid \alpha \leq \sigma )$ and for each $\alpha < \sigma$, $M_{\alpha + 1}/M_{\alpha} \cong C_\alpha$ for some $C_\alpha \in \mathcal C$. Since for each cardinal $\kappa$, all direct summands of $M_{\mathcal C}^{(\kappa)}$ are strong splitters, by induction on $\alpha < \sigma$, we get that $M_\alpha \cong \bigoplus_{\beta < \alpha} C_\beta$ is a strong splitter and $\Ext 1R{C_\alpha}{M_\alpha} = 0$, so the inclusion $M_\alpha \subseteq M_{\alpha + 1}$ splits. Hence $M = M_\sigma \cong \bigoplus_{\alpha < \sigma} C_\alpha$ is a strong splitter, too.   

Let $\kappa = \card R + \aleph_0$ and $\nu = \kappa + \card {M_{\mathcal C}}$.  

\begin{proposition}\label{A1}
The following conditions are equivalent:
	\begin{enumerate}
		\item[(1)] The class $\mathcal D$ is $\lambda$-categorical for all cardinals $\lambda > \nu$.
		\item[(2)] $\mathcal D$ is $\lambda$-categorical for some cardinal $\lambda \geq \nu$.
		\item[(3)] $M^{(\nu)} \cong N^{(\nu)}$ for all $0 \neq M,N \in \mathcal C$.
	\end{enumerate}
   \end{proposition}
\begin{proof} The proof is similar to the one for Proposition \ref{projSCC}, so we will only sketch it: Let $0 \neq M,N \in \mathcal C$. Then $\card M, \card N \leq \nu$. Assume (2) holds. Then  $M^{(\nu)}$ is isomorphic to a direct summand in $N^{(\nu)}$, and vice versa. Moreover, $M^{(\nu)} \cong N^{(\nu)} \oplus M^{(\nu)}$, and similarly for $N^{(\nu)}$. Hence $M^{(\nu)} \cong N^{(\nu)} \oplus M^{(\nu)} \cong N^{(\nu)}$ and (3) holds.

Assume (3). Let $D \in \mathcal D$ be of cardinality $\lambda > \nu$. By the above, $D$ is isomorphic to a direct sum of modules from $\mathcal C$, $D \cong \bigoplus_{C \in \mathcal C} C^{(\kappa_C)}$ for some cardinals $\kappa_C \geq 0$ ($C \in \mathcal C$). Let $D^\prime = \bigoplus_{C \in \mathcal C ^\prime} C^{(\kappa_C)}$, where $\mathcal C ^\prime$ is the set of all $C^\prime \in \mathcal C$ with $\kappa_C < \nu$. By (3), there is a module $0 \neq M \in \mathcal C$, such that for each $C \in \mathcal C \setminus \mathcal C ^\prime$, $C^{(\kappa_C)} \cong M^{(\kappa_C)}$. Since $\card D = \lambda > \card {\mathcal C} \geq \card {\mathcal C^\prime}$, we infer that $\card {D ^\prime} < \lambda$, and $D \cong D^\prime \oplus \bigoplus_{C \in \mathcal C \setminus \mathcal C ^\prime} M^{(\kappa_C)}$. As $\bigoplus_{C \in \mathcal C \setminus \mathcal C ^\prime} M^{(\kappa_C)}$ has cardinality $\lambda$, $D \cong D^\prime \oplus M^{(\lambda)}$. By (3), $D^\prime$ is isomorphic to a direct summand in $M^{(\lambda)}$, so by Eilenberg's Trick, $D \cong M^{(\lambda)}$, and (1) holds. 
\end{proof}    

In particular, if $\mu \geq \kappa$ is the least cardinal such that $\mathcal D$ is $\mu^+$-deconstructible then by the above either

\medskip
(i) $\Spec {\mathcal D} \subseteq \langle \mu, 2^\mu )$, and $\Spec {\mathcal D}$ is bounded, or else 

(ii) $\Spec {\mathcal D} \supseteq \langle (2^{\mu})^+, \infty )$, and for \emph{each} non-zero module $N \in \mathcal C$ and each cardinal $\lambda \geq (2^{\mu})^+$, all the elements of $\mathcal D$ of cardinality $\lambda$ are isomorphic to $N^{(\lambda)}$.    

Thus, in the Alternative 1, one can let {\lq}{\lq}sufficiently big{\rq}{\rq} mean {\lq}{\lq}of cardinality $\geq (2^{\mu})^+${\rq}{\rq}. 
 
\begin{remark}\label{remA1} The case (i) occurs, for example, when $\mathcal C = \{ C_1, C_2 \}$ where $C_1 \ncong C_2$ are countably generated indecomposable modules with local endomorphism rings and $M_C = C_1 \oplus C_2$ is a strong splitter. By Azumaya's theorem \cite[12.6]{AF}, $C_1^{(\lambda)} \ncong C_2^{(\lambda)}$ for any $\lambda \geq \kappa$, whence $\Spec {\mathcal D} = \emptyset$.

The case (ii) occurs, for example, when $\mathcal D = \Filt{\mathcal C}$ where $\mathcal C = \{ C \}$ and $C$ is a strong splitter. In this case $\mathcal D$ is the class of all modules isomorphic to $C^{(\rho)}$ for some cardinal $\rho \geq 0$, and $\kappa \leq \mu \leq \nu = \kappa + \card C$. Then $\Spec {\mathcal D} = \langle \nu, \infty )$ in case $C \cong C^{(\nu)}$, and $\Spec {\mathcal D} = \langle \nu ^+, \infty )$ when $C \ncong C^{(\nu)}$.
\end{remark}

\medskip
{\bf Alternative 2: $M_{\mathcal C}$ is not a strong splitter.} In this case, $\Spec {\mathcal D}$ is bounded by Lemma \ref{notail}. 

\medskip 
Now we arrive at the main result of this section:

\begin{theorem}\label{main} Let $\mathcal D$ be a deconstructible class of modules. Then the categoricity of $\mathcal D$ in a big cardinal is equivalent to its categoricity in a tail of cardinals. 
\end{theorem}
\begin{proof} Just observe that the equivalence holds in each of the Alternatives 1 and 2 discussed above: either $\Spec {\mathcal D}$ is bounded (in the case (i) of Alternative 1, and in Alternative 2) or else $\Spec {\mathcal D}$ contains a tail of cardinals (in the case (ii) of Alternative 1). 
\end{proof}

As pointed out by the referee, since for each cardinal $\mu$, there is only a set of $\mu^+$-deconstructible classes of modules over a fixed ring $R$, we can formulate Theorem \ref{main} in the following uniform way:

\begin{theorem}\label{main+} Let $R$ be a ring and $\mu$ be a cardinal such that $\mu \geq \card R + \aleph_0$. Then there exists a cardinal $\lambda_\mu$ such that for \emph{each} $\mu^+$-deconstructible class of modules $\mathcal D$, if $\mathcal D$ is $\lambda$-categorical for some $\lambda \geq \lambda_\mu$, then $\mathcal D$ is $\lambda$-categorical for every cardinal $\lambda \geq \lambda_\mu$.
\end{theorem}
 
\medskip
In the Alternative 1, we have seen that Theorem \ref{main+} holds for $\lambda_\mu = (2^{\mu})^+$. 

However, we do not have any a priori bound for $\lambda_\mu$ in the Alternative 2. Consider again the setting and notation of Lemma \ref{notail} -- so in particular, $M \in \mathcal D$ is not a strong splitter. Assume that $\Spec {\mathcal D}$ contains a cardinal $\lambda \geq \kappa$. Then $\lambda$ is unique. 

Using Lemma \ref{et} for $D = M$, we see that if $\lambda^{\kappa} = \lambda$, then both $M^{(\lambda)}$ and $N_\lambda$ are modules from $\mathcal D$ of cardinality $\lambda$, but $\Ext 1RM{N_{\lambda}} = 0$, while $\Ext 1RM{M^{(\lambda)}} \neq 0$, whence $M^{(\lambda)} \ncong N_\lambda$, a contradiction. In particular, the unique $\lambda$ cannot be of the form $2^\rho$ for any $\rho \geq \kappa$. Under GCH, we infer that $\lambda$ cannot be any successor cardinal $> \kappa$ (and in fact, any cardinal of cofinality $> \kappa$). 

So Lemma \ref{et} gives a number of restrictions for $\lambda$. Moreover, in most particular cases, $\Spec {\mathcal D}$ will not contain any cardinal $\lambda \geq \kappa$ (e.g., when $\mathcal D$ contains a module $N$ of cardinality $\leq \kappa$ such that $N \notin \Add M$).
   
We will finish by showing that a simple modification of the proof of Lemma \ref{notail} yields yet another restriction for $\lambda$ in ZFC: if $\lambda > \kappa$, then $\lambda$ must be a singular cardinal:

\begin{theorem}\label{notregular} Let $R$ be a ring, and $\mathcal D$ be a deconstructible class of modules containing a module $M$ which is not a strong splitter. Let $\kappa = \card M + \card R + \aleph_0$. Then $\mathcal D$ is not $\lambda$-categorical for \emph{any} regular cardinal $\lambda > \kappa$. 
\end{theorem}
\begin{proof} By Lemma \ref{single}, $\Ext 1RM{M^{(\nu)}} \neq 0$ for each cardinal $\nu \geq \kappa$. Let $\lambda > \kappa$ be a regular cardinal. By induction, we define a strictly increasing chain of modules $( M_\alpha \mid \alpha < \lambda )$ of cardinality $< \lambda$ as follows: $M_0 = M^{(\kappa)}$; if $M_\alpha$ is defined and $M \cong M^{(\nu)}$ for some $\kappa \leq \nu < \lambda$, then we define $M_{\alpha + 1}$ using a non-split extension $0 \to M_{\alpha} \to M_{\alpha +1} \to M \to 0$ which exists because $\Ext 1RM{M^{(\nu)}} \neq 0$. If $M_\alpha$ is defined and $M \ncong M^{(\nu)}$ for any $\kappa \leq \nu < \lambda$, we let $M_{\alpha +1} = M_{\alpha} \oplus M$. If $\alpha < \lambda$ is a limit ordinal, we put $M_\alpha = \bigcup_{\beta < \alpha} M_\beta$. Let $N = \bigcup_{\alpha < \lambda} M_\alpha$. Then $N \in \mathcal D$ and $\card N = \lambda$.

Assume $\mathcal D$ is $\lambda$-categorical. Then $N \cong M^{(\lambda)}$, so there exists an $R$-independent family $( P_\alpha \mid \alpha < \lambda)$ of submodules of $N$ such that $N = \bigoplus_{\alpha < \lambda} P_\alpha$ and $P_\alpha \cong M$ for each $\alpha < \lambda$. 

By induction, we define an increasing chain of ordinals $( \beta_n \mid n < \omega )$ and an increasing chain of subsets $( I_n \mid n < \omega)$ of $\lambda$ of cardinality $< \lambda$ as follows: $\beta_0 = 0$ and $I_0$ is the least subset of $\lambda$ such that $M_{\beta_0} \subseteq \bigoplus_{\alpha \in I_0} P_\alpha$. Since $M_{\beta_0} = M^{(\kappa)}$, $\card (I_0) \leq \kappa$. For the inductive step, if $\card (I_n) < \lambda$, then the regularity of $\lambda$ implies existence of an ordinal $\beta_n < \beta_{n+1} < \lambda$ such that $\bigoplus_{\alpha \in I_n} P_\alpha \subseteq M_{\beta_{n+1}}$. Since $\card M_{\beta_{n+1}} < \lambda$, there is a subset $I_n \subseteq I_{n+1} \subset \lambda$ such that $\card (I_{n+1}) < \lambda$ and $M_{\beta_{n+1}} \subseteq \bigoplus_{\alpha \in I_{n+1}} P_\alpha$.    

Let $\beta = \sup_{n < \omega} \beta_n < \lambda$, and $I = \bigcup_{n < \omega} I_n$. Then 
$$M_\beta = \bigcup_{n < \omega} M_{\beta_n} = \bigcup_{n < \omega} (\bigoplus_{\alpha \in I_n} P_\alpha) = \bigoplus_{\alpha \in I} P_\alpha.$$ 
So $M_\beta$ is a direct summand in $N$, and hence in $M_{\beta + 1}$. However, $M_\beta \cong \bigoplus_{\alpha \in I} P_\alpha$ and $\lambda > \card I \geq \card I_0 \geq \kappa$, so $M_\beta \cong M^{(\nu)}$ for some $\kappa \leq \nu < \lambda$. Hence $M_\beta$ is not a direct summand in $M_{\beta +1}$ by the inductive step of our construction of the module $N$, a contradiction.   
\end{proof}

The following questions remain open:

\medskip
1. Is there an instance of Lemma \ref{notail} where $\card \hbox{Spec}_{\kappa}(\mathcal D) = 1$? If so, then $\hbox{Spec}_{\kappa}(\mathcal D) = \{ \lambda \}$ for a cardinal $\lambda \geq \kappa$, and $\lambda$ is singular in case $\lambda > \kappa$ by Theorem \ref{notregular}.

2. If so, how does then the value of $\lambda$ depend on $\mu$ (= the least infinite cardinal $\geq \card R$ such that $\mathcal D$ is $\mu^+$-deconstructible)? The answer would help to determine the value of $\lambda_\mu$ in Theorem \ref{main+}. 

\medskip
\begin{acknowledgment}
The author thanks Kate\v{r}ina Fukov\'{a} for valuable comments on an earlier draft of this paper, and the referee for several suggestions that helped to improve the paper.    
\end{acknowledgment}

\end{document}